\documentclass[12pt]{amsart}

\usepackage{graphicx}
\usepackage{psfrag}

\newcommand{\Lat}{\mathcal{L}}
\newcommand{\Int}{\mathbf{Z}}
\newcommand{\nint}{\mathbf{N}}
\newcommand{\vrtx}{V_\Lat}
\newcommand{\prd}{P_\Lat}
\newcommand{\inI}{in_\omega(I_\Lat)}
\newcommand{\lcm}{{\rm lcm}}
\newcommand{\triang}{\Delta_\omega}
\newcommand{\Vrt}{ {\rm Vert}}
\newcommand{\conv}{ {\rm conv}}

\newtheorem{theorem}{Theorem}[section]
\newtheorem{lemma}[theorem]{Lemma}

\newtheorem{proposition}[theorem]{Proposition}

\newtheorem{cor}[theorem]{Corollary}

\theoremstyle{definition}

\newtheorem{definition}[theorem]{Definition}
\newtheorem{example}[theorem]{Example}
\newtheorem{remark}[theorem]{Remark}

\def\bartau{{\bar \tau}}

\begin{document}
\title{The Vertex Ideal of a Lattice }
\author{ Serkan Ho\c{s}ten}

\address{Department of Mathematics, San Francisco State University, San Francisco, CA 94132}

\email{ serkan@math.sfsu.edu}

\author{Diane Maclagan}
\address{ School of Mathematics, Institute for Advanced Study, Princeton, NJ 08540}

\email{maclagan@ias.edu}

\begin{abstract}
We introduce a monomial ideal whose standard monomials
encode the vertices of all fibers of a lattice. We study the minimal 
generators, the radical, the associated primes and the primary decomposition 
of this ideal, as well as its relation to initial ideals of lattice ideals.
\end{abstract}

\maketitle

\section{Introduction}

The main purpose of this paper is to introduce and study a monomial ideal, the
vertex ideal, associated to a sublattice of $\Int^n$.  We relate
algebraic properties of this ideal to combinatorial properties of the
lattice.

\begin{definition}
Let $\Lat$ be a lattice in $\Int^n$ with $\mathrm{dim}(\Lat) = m$.
For $u \in \nint^n$ we define $P_u := \conv\{v \in \nint^n \,: \,\, u-v \in 
\Lat \}$ to be the {\em fiber} of $u$ with respect to $\Lat$. Clearly, if
$v \in P_u$ then $P_u = P_v$.  
\end{definition}

Each fiber $P_u$ is a rational polyhedron, by Theorem 16.1 in
\cite{Sch}, and hence has finitely many vertices $\Vrt(P_u)$. We start
with the observation (Proposition \ref{vldefn}) that the union of all
$\Vrt(P_u)$, $u \in \nint^n$ forms an order ideal of $\nint^n$.  We
call the monomial ideal which is the complement of this order
ideal the {\em vertex ideal} of $\Lat$, and denote it by  $\vrtx$.

One motivation for studying vertex ideals comes from the theory of
integer programming.  Suppose $A \in \nint^{d \times n}$ is a matrix
of rank $d$ with no zero columns.  Let $\nint A$ be the submonoid of
$\nint^d$ consisting of nonnegative integer combinations of the
columns of $A := [a_1\ldots a_n]$.  Integer programming is concerned
with minimizing a fixed linear form $c \cdot x$, where $c \in
\mathbf{R}^n$, over $\{u \in \nint^n: \,\, Au = b\}$ for a fixed $b
\in \nint A$.  Note that if we let $\Lat = ker(A) \cap \Int^n$, then
for $v \in \nint^n$, the fiber of $v$ is a polytope $\conv\{u \in
\nint: \,\, Au = Av\}$ (in this case we denote the fiber of $v$ by 
$P_b$ where $b = Av$). Hence studying the vertex ideal $\vrtx$ in
this context gives information about the vertices of all integer
programming polytopes as $b$ varies in $\nint A$.  Commutative algebra
and computational algebraic geometry enter this picture through the
connection between integer programming and Gr\"obner bases and initial
ideals of the {\it toric ideal} of $A$ (see \cite{GBCP}, \cite{ST} and
\cite{Rek}).

A second motivation comes from the recent work of Saito, Sturmfels and 
Takayama \cite{SST} on hypergeometric differential equations. One observation
these authors make is that the set of all {\em generic $A$-hypergeometric 
series solutions} to a GKZ $A$-hypergeometric system is indexed by the 
top-dimensional {\em standard pairs} of $\vrtx$ where $\Lat = ker(A) 
\cap \Int^n$ (pp. 129-131 in \cite{SST}).
This leads us to studying the standard pairs (and hence the associated
primes) of $\vrtx$. 
  
In Section 2 of this paper we start by giving a naive algorithm to
construct $\vrtx$ in Theorem \ref{vertex}.  This first algorithm needs
{\em all} initial ideals of the associated {\it lattice ideal}
$I_\Lat$, and therefore it is highly inefficient for large problems.
We remedy this by giving an improved algorithm to construct a
generating set for $\vrtx$, using the {\em Graver basis} elements of
$I_\Lat$.  This second algorithm depends on a characterization of
$\vrtx$ which is derived from only the geometric properties of the
lattice.  We also describe the radical of $\vrtx$ as the
{\em Stanley-Reisner ideal} of a matroid complex.

In Section 3 the second motivation we cited above for studying $\vrtx$
leads us to investigate the associated primes of $\vrtx$. First we
give a characterization of $Ass(\vrtx)$ and compute the irreducible
primary decomposition of $\vrtx$ in terms of a family of polytopes
using similar methods to those found in \cite{HT}. This allows us to
give some necessary conditions for a prime being an embedded
associated prime of $\vrtx$ when $\mathrm{dim}(\Lat) = 2$. In
particular, we show that the irrelevant maximal ideal $\langle x_1,
\ldots, x_n \rangle$ cannot be associated to $\vrtx$ in this
case. This result fails when $\mathrm{dim}(\Lat) \geq 3$, and we give
a counterexample. This seemingly harmless counterexample turns out to
be a very interesting one for our first motivation, integer
programming.  It provides a counterexample to a conjecture about the
complexity of codimension three integer programs. More precisely, it
gives a counterexample to Conjecture 6.1 in \cite{ST} which
hypothesized that every cone in the {\em Gr\"obner fan} of a
codimension three toric ideal has at most four facets.

In Section 4 we  define another monomial ideal, $\prd$, closely related to 
$\vrtx$. We show that the {\em product ideal} $\prd$ has the same 
radical as $\vrtx$. In two interesting special cases, we prove that 
$\prd$ (which is easier to compute) is equal to $\vrtx$. The first case is
when $\Lat$ comes from  a unimodular matrix $A$. The second case is
when $\Lat \subseteq \Int^2$ and $\mathrm{dim}(\Lat) = 2$. This implies that for any
two dimensional lattice, we have  $Top(\prd) = Top(\vrtx)$.

\section{The Minimal Generators and the Radical of the Vertex Ideal}

The first goal of this section is to come up with useful characterizations of 
$\vrtx$ which we  use for devising a relatively efficient algorithm.
We then give a combinatorial description of the radical of the
vertex ideal. We first show the existence of the vertex ideal.  

\begin{proposition} \label{vldefn}
Let $\Lat$ be a lattice in $\Int^n$, and let $P_u$ be a fiber of $\Lat$.  
For any vertex $v$ of $P_u$, if $v_i > 0$, 
then $v-e_i$ is a vertex of $P_{u-e_i}$ where $e_i$ is the $i$-th unit 
vector. In other words, there exists a monomial ideal $\vrtx$ in 
$S=k[x_1,\ldots,x_n]$  where $x^v \notin \vrtx$ if and only if 
$v \in \Vrt(P_u)$ for a fiber $P_u$ of $\Lat$.
\end{proposition}
\begin{proof} If $v-e_i$ is not a vertex of $P_{u-e_i}$, then it is 
in the convex hull of vertices $v_1', v_2', \ldots, v_k'$ of $P_{u-e_i}$.
But then $v$ would be in the convex hull of 
$v_1'+e_i, v_2'+e_i, \ldots, v_k'+e_i$. This contradiction proves the
first statement, and hence implies that the union of all
$\Vrt(P_u)$,  $u \in \nint^n$ forms an order ideal of $\nint^n$. This
is equivalent to the second statement.
\end{proof}

We now give a first algorithm to compute $V_\Lat$.  To do this, we
first associate a binomial ideal to $\Lat$.

\begin{definition}
The lattice ideal $I_\Lat$ is defined by 
$$I_\Lat = \langle x^u-x^v \,\,: \,\, u,v \in \nint^n, \,\, u-v \in \Lat 
\rangle.$$
\end{definition}

Lattice ideals have been widely studied, see for example 
\cite{HT}, \cite{PStu}, \cite{PStu2}.  In this context we are interested in the
initial ideals of $I_{\Lat}$. For a {\em weight vector} 
$\omega \in \mathbf{R}^n$ such that $\omega \cdot u > 0$ 
for every non-zero vector $u \in \nint^n \cap \Lat$, we let $\inI$ be
the ideal $\langle in_\omega(f) \,\,: \,\, f \in I_\Lat \rangle$ where
$in_\omega(f)$ is the sum of all terms of $f$ with maximum
$\omega$-value. If the {\em initial ideal} $\inI$ is a monomial ideal
we call $\omega$ a {\em generic} weight vector.  Our assumption on 
$\omega$ ensures that each fiber $P_u$ has a bounded face which
minimizes the linear functional $\omega \cdot x$. Then the genericity of
$\omega$ is equivalent to the condition that each such bounded face 
is a vertex $v$ of $P_u$.  

\begin{theorem}\label{vertex} The vertex ideal $\vrtx$ is equal to
$\bigcap_{\omega} \inI$ where $\omega$ is a generic weight vector. 
\end{theorem}

\begin{proof}Since for any
two lattice points $u, v \in P_u$ we have $u-v \in \Lat$,  a
monomial is a standard monomial of $\inI$ if and only if its exponent vector
minimizes the linear functional $\omega \cdot x$ in $P_u$ \cite{SWZ}. Hence 
the monomial $x^v$ is a standard monomial of
$\bigcap_{\omega} \inI$ if and only if $v$ is the minimizer of 
$\omega \cdot u$ for $u \in P_v$ for {\em some} generic weight vector. But 
these are precisely the vertices of the fibers of $\Lat$. 
\end{proof}

Using this theorem we have a first algorithm for computing $\vrtx$:
compute all initial monomial ideals of $I_\Lat$ and take their
intersection.  We note that this is a finite algorithm, as any ideal in $S$
has only a finite number of different initial ideals.  The list of all
initial ideals of $I_{\Lat}$ can be computed with the software {\tt
TiGERS} \cite{HuT}.  This first algorithm is not, however, completely
satisfactory, as the number of initial ideals can be much larger than
the subset needed to define the intersection.  In order to illustrate
this point we use the following example, where the number of initial
ideals depends exponentially on the data of the lattice.

\begin{example} \label{2x2} Let $I_\Lat$ be the ideal generated by the 
$2 \times 2$ minors of a generic $2 \times n$ matrix $X = (x_{ij})$.
This is a prime lattice ideal which is the defining ideal of the Segre
embedding of $\mathbf{P}^1 \times \mathbf{P}^{n-1}$ into
$\mathbf{P}^{2n-1}$. Proposition 5.4 of \cite{GBCP} shows that with
respect to the reverse lexicographic term order $x_{11} \prec \cdots
\prec x_{1n} \prec x_{21} \prec \cdots \prec x_{2n}$, these $2 \times
2$ minors form a reduced Gr\"obner basis. By permuting columns of $X$,
and using the corresponding reverse lexicographic term order, one gets
$n!$ distinct initial ideals.  This shows that $I_\Lat$ has at least
$n!$ initial ideals. In Remark \ref{segreegcont} we will see
that the vertex ideal can be constructed as the intersection of only
$n2^{n-1}$ initial ideals.  As $\frac{n2^{n-1}}{n!} \rightarrow 0$ as
$n \rightarrow \infty$, a vanishingly small proportion of the initial
ideals are needed to construct $V_{\Lat}$ in this family.

\end{example}

  Below we give a more efficient description of the minimal generators
of the vertex ideal.  For this, we need to define the Graver basis of
$\Lat$.

\begin{definition} \label{graver}
Suppose $\Lat \subseteq \Int^n$ and let $\mathbf{R}_\rho$ be the orthant
defined by the sign pattern $\rho \in \{+,-\}^n$. Then 
$\Lat \cap \mathbf{R}_\rho$ is a finitely generated monoid with 
a unique minimal generating set $H_\rho$, its {\em Hilbert basis} 
(Theorem 16.4 in \cite{Sch}). The {\em Graver basis} $Gr_\Lat$ 
of $\Lat$  (or $I_\Lat$) is defined to be the union of all such $H_\rho$.  
\end{definition}

\begin{lemma}\label{circuits}
Suppose $\sum_i c_i(\alpha_i-\beta_i)=0$, where $c_i>0$, and
$\alpha_i-\beta_i \in Gr_\Lat$ with $\alpha_i, \beta_i \in \nint^n$ and
$supp(\alpha_i) \cap supp(\beta_i) = \emptyset$.  Then $x^v = 
\lcm_i(x^{\alpha_i})$ is in $\vrtx$.
\end{lemma}
\begin{proof}
Suppose $x^v$ is not in $\vrtx$.  This means that $v$ is a vertex of
$P_v$, so there is some $\omega \in \mathbf R^n$ such that $\omega
\cdot v > \omega \cdot u$ for all lattice points $u \in P_v \setminus
\{v\}$.  But now 
$$\omega \cdot v > \omega \cdot (v- (\alpha_i-\beta_i)) \mbox{   for each } i$$
because $\alpha_i \leq v$ means that
$v^{\prime}=(v-\alpha_i+\beta_i) \in \mathbf N^n$, and thus
$v^{\prime}$ is a lattice point in $P_u \setminus \{v\}$.  This
implies

\begin{eqnarray*}
\sum_i \omega \cdot (c_iv) & > & \sum_i \omega \cdot (c_iv-c_i(\alpha_i-\beta_i))\\
& = &  \sum_i \omega \cdot (c_iv) - 
\omega \cdot \sum_i c_i(\alpha_i-\beta_i)  \\
& = & \sum_i \omega \cdot (c_iv)\\
\end{eqnarray*}
This contradiction shows that $x^v$ is in $\vrtx$.
\end{proof}

\begin{cor}\label{mingens} The minimal generators of $\vrtx$ are
of the form as in Lemma \ref{circuits}.
\end{cor}

\begin{proof} Let $x^u$ be a minimal generator of $\vrtx$. Hence 
$u$ is not a vertex of its fiber, and therefore it is a convex
combination $\sum_i\lambda_i v_i$ of some vertices $v_i$ of $P_u$,
where $0 \leq \lambda_i \leq 1$ and $ \sum_i \lambda_i=1$.  Since
$u-v_i$ is in $\Lat$, we have $u - v_i = \sum_j c_{ij} (\alpha_{ij} -
\beta_{ij})$ where $\alpha_{ij} - \beta_{ij}$ are Graver basis
elements with $\alpha_{ij} \leq u$ and $\beta_{ij} \leq v_i$, and
$c_{ij} \in \mathbf Z_{\geq 0}$. Now clearly $\sum_i \lambda_i(u-v_i) = 0$, and
thus $\sum_{i,j} \lambda_i c_{ij}(\alpha_{ij} - \beta_{ij}) = 0$.  By
Lemma \ref{circuits}, $x^v = \lcm_{ij}(x^{\alpha_{ij}})$ is in $\vrtx$.  But
$x^v$ divides $x^u$, and $x^u$ is a minimal generator, so $x^u = x^v$.
\end{proof}

Corollary \ref{mingens} implies that the minimal generators of $\vrtx$
can be computed by identifying all positive linear dependencies
among Graver basis elements of $\Lat$. In fact only the
minimal positive dependencies, known as {\em positive circuits}, are
needed. We summarize this as follows.

\begin{theorem} \label{pos-circuits}
Let $Gr_\Lat= \{\alpha_i-\beta_i\} $ be an {\em ordered} Graver basis
of $\Lat$, so that $\alpha-\beta \in Gr_\Lat$ implies $\beta-\alpha
\in Gr_\Lat$.  If $\tau$ is the support of a positive circuit 
$\sum_{i \in \tau} c_i(\alpha_i-\beta_i)=0$
we define
$x^{m_{\tau}}$ to be $\lcm_{j \in \tau} x^{\alpha_j}$.  Then
$$\vrtx= \langle x^{m_{\tau}} | \tau \mbox{ is the support 
of a positive circuit of } Gr_\Lat \rangle.$$
\end{theorem}

\begin{proof}
If $\tau$ is the support of a positive circuit of $Gr_\Lat$,
Lemma \ref{circuits} implies that $x^{m_{\tau}}$ is in $\vrtx$. And Corollary \ref{mingens} says that
every minimal generator of $\vrtx$ is of this form.
\end{proof}

Theorem \ref{pos-circuits} gives our second, more efficient, algorithm 
to compute $\vrtx$: after computing the Graver basis $Gr_\Lat$, identify 
each positive circuit $\tau$ of $Gr_\Lat$  and compute 
$x^{m_{\tau}}=\lcm_{j\in \tau} x^{\alpha_j}$.   

We observe that not all vectors of $Gr_\Lat$ are necessary.  When
$\Lat \cap \nint^n = \{0\}$, it suffices
to replace $Gr_\Lat$ by the ordered {\em universal Gr\"obner basis} of
$\Lat$.  See \cite[Chapter 7]{GBCP} for information on computing the universal
Gr\"obner basis.

The next result in this section describes the radical of $\vrtx$. Let
$B \in \Int^{n \times m}$ be a matrix whose columns form a basis for
the $m$-dimensional lattice $\Lat$. We will denote the rows of $B$ by
$b_1, \ldots, b_n$. Now if $\omega$ is a generic cost vector, the
vector $\omega B$ is contained in the relative interior of a set
$\mathcal{C}$ of $m$-dimensional simplical cones with generators from
$\{b_1, \ldots, b_n\}$. We define $\triang$ to be the simplicial
complex generated by the complementary indices of the generators of
the cones in $\mathcal{C}$. By its definition, $\triang$ is an
$(n-m)$-dimensional pure simplicial complex on $\{1,\ldots, n\}$. We
also note that this simplicial complex is the {\em regular
triangulation} of $A$ with respect to $\omega$ when $\Lat = ker(A)
\cap \Int^n$ (see Chapter 8 in \cite{GBCP}). Extending the connection
between Stanley-Reisner ideals of regular triangulations of $A$ and
the radicals of the initial ideals of $I_\Lat$, we get the following
proposition (Corollary 2.9 in \cite{HT}, see also Section 7 in
\cite{SWZ}).  Recall that the Stanley-Reisner ideal of a simplicial
complex is the ideal generated by the minimal non-faces of the
complex.

\begin{proposition}\label{St-Reis} The radical of $\inI$ is the 
Stanley-Reisner ideal of the simplicial complex $\triang$.
\end{proposition}

Now we are ready to prove the following theorem:
\begin{theorem}\label{radva}
The radical of $V_\Lat$ is $\bigcap_{\sigma} \langle x_i
: i \in \sigma \rangle $ where the intersection is over all 
linearly independent subsets of $\{b_1, \ldots, b_n\}$ of size $m$.
\end{theorem}
\begin{proof}
\begin{eqnarray*}
{\mathrm rad}(V_\Lat)&=&{\mathrm rad}(\bigcap_{\omega \ generic} \inI) \\
	&=&\bigcap_{\omega \ generic} {\mathrm rad}
		(\inI) = \bigcap_{\triang} I_{\triang} \\
	&=&\bigcap_{\triang} \bigcap_{\tau
		\in \triang} \langle x_i : i \notin \tau \rangle \\
	&=&\bigcap_{\sigma: \mathrm{dim}(\sigma) = m} 
\langle x_i : i \in \sigma \rangle \\
\end{eqnarray*}
where $I_{\triang}$ is the Stanley-Reisner ideal of $\triang$.  
We have the first equality on the second line because
taking the radical commutes with intersections, while the second equality follows
from Proposition \ref{St-Reis}. The third line is a standard
result on Stanley-Reisner ideals, and the last line follows because
the complement of the indices of the generators of any full dimensional 
simplicial cone $\{b_{i_1}, \ldots, b_{i_m} \}$ is involved in some $\triang$.
\end{proof}

This result can be interpreted using the notion of a matroid complex.

\begin{definition}
The matroid complex $\Delta(\mathcal{M})$ of a matroid $\mathcal{M}$ 
is the simplicial complex where the simplices are the independent sets of 
$\mathcal{M}$.
\end{definition}

If $\Lat \subset \Int^n$ is a lattice of dimension $m$ generated by
the columns of a matrix $B \in \Int^{n \times m}$, then the
complements of {\em bases} (i.e. linearly independent subsets of rows
of size $m$) of $B$ form the maximal independent sets of a matroid
$\mathcal{M}(\Lat)$.  Hence the matroid complex
$\Delta(\mathcal{M}(\Lat))$ is the simplicial complex whose maximal
simplices are the union of the maximal simplices occuring in
$\Delta_\omega$ for all generic $\omega$. Note that when $\Lat =
ker(A) \cap \Int^n$ for a matrix $A$, then $\mathcal{M}(\Lat)$ is the
matroid of all linearly independent subsets of the columns of $A$, and
$\Delta(\mathcal{M}(\Lat))$ is the simplicial complex whose maximal
simplices are the union of all maximal simplices appearing in the
regular triangulations of $A$. We now get the following corollary.

\begin{cor}\label{matroid}
The Stanley-Reisner ideal of $\Delta(\mathcal{M}(\Lat))$ is the radical of 
$V_\Lat$.
\end{cor}

\begin{proof}
The Stanley-Reisner ideal of $\Delta(\mathcal{M}(\Lat))$ is
$$I_{\Delta(\mathcal{M}(\Lat))} = \bigcap_{\tau \in \Delta(\mathcal{M}(\Lat))} 
\langle x_i : i \notin \tau \rangle. $$
The above intersection can be taken over all $\tau$ where $\tau$ is a maximal
face. Then since $\tau \in \Delta(\mathcal{M}(\Lat))$ if and only 
if $\{b_i \,: i \notin \tau\}$ forms a basis of $B$ where $B$ is matrix 
whose columns are a basis for $\Lat$ , Theorem \ref{radva} implies that 
$I_{\Delta(\mathcal{M}(\Lat))} = V_\Lat$.
\end{proof}

\begin{remark} \label{segreegcont}
 We can now prove the last claim in Example \ref{2x2}.  The vertex
ideal $\vrtx$ of the $2 \times 2$ minors of a generic $2 \times n$
matrix is a radical ideal, as all the initial ideals are radical,
because the corresponding configuration is unimodular.  Hence we can
use the intersection formula in the proof above. The maximal faces
over which we need to take the intersection are determined by maximal
independent sets of the collection $\{ e_i \oplus e_j : i=1,2, \mbox{
and } 1 \leq j \leq n \}$.  These are in bijection with the distinct
spanning trees of the complete bipartite graph $K_{2,n}$.  There are
$n2^{n-1}$ such spanning trees, as exactly one vertex in the $n$-block
is connected to both vertices in the 2-block.
\end{remark}

Finally, we observe that the Hilbert series of $\vrtx$ gives us
information about the number of vertices of the fibers $P_u$ of
$\Lat$.

\begin{proposition} \label{Hilbert}The Hilbert series $H(S/\vrtx; z_1,\ldots,z_n)$ of 
$S/\vrtx$ is $\sum_u {\bf z}^u$,
where the sum is taken over all vertices $u$ of all fibers $P_u$. 
When $\Lat = ker(A) \cap \Int^n$ for an integer matrix $A = [a_1,\ldots,a_n]$ 
then
$$H(S/\vrtx; z_1{\bf t}^{a_1},\ldots, z_n{\bf t}^{a_n}) =
\sum_{b \in \nint A} (\sum_{u \in \mathrm{Vert}(P_b)} {\bf z}^u) \cdot 
{\bf t}^b, \,\,\, and,$$
$$H(S/\vrtx; {\bf t}) = \sum_{b \in \nint A} 
|\mathrm{Vert}(P_b)| \cdot {\bf t}^b.$$
\end{proposition}

We can derive information about the fibers $P_u$ from the Hilbert
function for $\vrtx$.  An example is given in the following
proposition.

\begin{proposition}
If $\Lat=ker(A) \cap \mathbf Z^n$ for a $1 \times n$ matrix 
$A = [a_1, \ldots, a_n]$ where $a_i \in \nint$, then
the number of vertices of a fiber $P_u$ is eventually periodic, with
period dividing $\lcm_i(a_i)$.

\end{proposition}

\begin{proof}
The Hilbert series $H(S/\vrtx;t)$ can be written in the form 
$$\frac{p(t)}{\prod_{i=1}^n (1-t^{a_i})},$$ for some polynomial $p(t)$.
This means that the Hilbert function of $S/\vrtx$ at $b$, which counts
the number of vertices of $P_u$ when $Au=b$, eventually agrees with a
quasi-polynomial evaluated at $b$.  As there is an upper bound, given
by the number of initial ideals of $I_{\Lat}$, on the number of
vertices of any $P_u$, this polynomial part of the quasi-polynomial
must be a constant.  As the period of the quasi-polynomial divides $\lcm(a_i)$, the result follows.

\end{proof}

We observe that a more constructive proof of this proposition can also
be given using the notion of {\em atomic fibers}, defined in
\cite{AHLM}.
 
\section{Associated Primes and Standard Pairs of $\vrtx$}

With the relation between initial ideals and $V_{\Lat}$ given in
Theorem \ref{vertex} it is natural to ask which properties of the
initial ideals of a lattice ideal pass to $\vrtx$. For example, these
initial ideals possess the rare property that their associated primes
come in saturated chains \cite{HT}.  Although we do not determine if
this property holds for the vertex ideal, this section provides some
tools for approaching this question.  Furthermore, while investigating
the associated primes of $\vrtx$, we construct a lattice which
provides a counterexample to a conjecture about codimension three
toric ideals.

Since $\vrtx$ is a monomial ideal, all of its associated primes are
monomial primes of the form $\mathcal{P}_{\sigma}= \langle x_i : i
\not \in \sigma \rangle$ where $\sigma \subseteq [n]
:=\{1,\ldots,n\}$. 

\begin{lemma}\label{associatedintersect}
The set of associated primes $\mathrm{Ass}(\vrtx)$ of $\vrtx$ is 
contained in $\cup_{\omega} \mathrm{Ass}(\inI)$, the union of the associated
primes of all initial ideals of $I_\Lat$. Furthermore, the set of 
minimal primes of $\vrtx$ is precisely the union of the minimal primes of 
all initial ideals of $I_\Lat$. 
\end{lemma}

\begin{proof} Using Theorem \ref{vertex}, the first statement follows from 
the fact that if two ideals $I$ and $J$ have minimal primary
decompositions $\cap_i \mathcal{P}_i$ and $\cap_j \mathcal{P}_j'$,
then $(\cap_i \mathcal{P}_i) \cap (\cap_j \mathcal{P}_j')$ is a (not
necessarily minimal) primary decomposition of $I \cap J$.  Minimal
primes of a intersection of monomials ideals are always contained in
the union of the minomial primes of the ideals.  The fact that this
containment is an equality in this case follows from the fact, used in
Theorem \ref{radva}, that all minimal primes of all initial ideals
have the same dimension.
\end{proof}

\begin{example} The associated primes of $\vrtx$ can be strictly contained
in $\bigcup_{\omega} \mathrm{Ass}(\inI)$.  Consider the matrix $A=[1
\, 2 \, 3]$ and $\Lat = ker(A) \cap \mathbf Z^n$.  For this lattice,
$\vrtx= \langle abc,a^2b,a^3c,b^3c^2 \rangle$, which has primary
decomposition $\langle a^3,ab,b^3 \rangle \cap \langle
a^2,ac,c^2\rangle \cap \langle b,c \rangle$, so the associated primes
of $\vrtx$ are $\langle a,b \rangle$, $\langle a,c \rangle$ and
$\langle b,c \rangle$.  For $\omega=(100,10,1)$ $\inI=\langle a^2, ab,
ac, b^3 \rangle$.  This has primary decomposition $\langle a , b^3
\rangle \cap \langle a^2,b,c \rangle$, so we have $\langle a,b,c
\rangle \in \cup_{\omega} \mathrm{Ass}(\inI)$.
\end{example}

Corollary 3.5 of \cite{HT} gives bounds on the dimensions and
codimensions of initial ideals of $I_\Lat$.  Combined with Lemma
\ref{associatedintersect} we get the following fact about the dimension
and codimension of the associated primes of $\vrtx$.

\begin{proposition}\label{dimension}
The dimension of an associated prime of $\vrtx$ for a lattice of
dimension $m$ is at least $\max(0,n-(2^m-1))$ and the codimension is
at most $\min(n,2^m-1)$. 
\end{proposition}

For our purposes it is more convenient to study the associated 
primes of $\vrtx$
via its {\em standard pairs} \cite{STV}.  For a vector $u \in \mathbf
N^n$ we denote by $\mathrm{supp}(u)$ the set $\{ i : u_i \neq 0 \}$.

\begin{definition}
An {\em admissible pair} of a monomial ideal $M$ is a pair $(x^u,\tau)$ with 
$\tau \subseteq [n]$ such that $\mathrm{supp}(u) \cap \tau = \emptyset$, 
and $x^{u+v} \not \in M$ for all $v$ 
with $\mathrm{supp}(v) \subseteq \tau$.  We place a partial order on the 
set of admissible pairs of $M$ by declaring $(x^u,\tau) \prec (x^v,\sigma)$ if 
$x^v | x^u$ and $\mathrm{supp}(u-v) \cup \tau \subseteq \sigma$.  
The maximal elements of the set of admissible pairs with respect to 
this order are called {\em standard pairs}.
\end{definition} 

In the rest of the paper we  use a polyhedral characterization of the 
standard pairs of $\vrtx$ following the results and terminology 
in \cite{HT} and
\cite{HT2}. We start with a characterization which follows from the
definition of standard pairs.

\begin{proposition} \label{stdpair1}
The pair $(x^u,\tau)$ is a standard pair of $\vrtx$ if and only if $u$
is a vertex of $P_u$, $\mathrm{supp}(u) \cap \tau = \emptyset$,
$u+v$ is a vertex of $P_{u+v}$ for all $v$ with $\mathrm{supp}(v)
\subseteq \tau$, and for all $i \not \in \tau$ there is some
$v^{\prime}$ with support in $\tau \cup \{i\}$ such that
$u+v^{\prime}$ is not a vertex of $P_{u+v^{\prime}}$.
\end{proposition}

As in the previous section, let $B \in \Int^{n \times m}$ such that
the columns of $B$ form a lattice basis for $\Lat$. Given $u \in
\nint^n$, we can define the polyhedron $Q_u := \{x \in \mathbf{R}^m:
Bx \leq u\}$.  The lattice points in $Q_u$ and the lattice points in
$P_u$ are in bijection by the correspondence $z \in Q_u \cap \Int^m
\longleftrightarrow u-Bz \in P_u$. The origin of $\Int^m$ is in $Q_u$
for all $u \in \mathbf N^n$ and corresponds to $u \in P_u$. We let
$R_u$ be the convex hull of the lattice points in $Q_u$.  Note that
$R_u$ is affinely isomorphic to $P_u$.  For a subset $\tau \subseteq
[n]$ we denote by $\bartau$ the complement of $\tau$, so $\bartau=[n]
\setminus \tau$.  With this convention we define $Q^{\bartau}_u$ to be
the polyhedron $\{x \in \mathbf{R}^m: B^{\bartau}x \leq u^{\bartau}\}$
where the inequalities defining $Q_u$ corresponding to $\tau$ are
omitted.  $R^{\bartau}_u$ denotes the convex hull of the lattice
points in $Q^{\bartau}_u$.  We now reformulate the characterization of
standard monomials and standard pairs of $\vrtx$.

\begin{theorem} \label{stdpair2}
The monomial $x^u$ is a standard monomial of $\vrtx$ if and only
if the origin is a vertex of $R_u$. Moreover, a pair $(x^u, \tau)$ 
is a standard pair of $\vrtx$ if and only if the origin is a vertex of 
$R^{\bartau}_u$ and it is {\bf not}
a vertex of $R^{\bartau \setminus i}_u$ for any $i \in \bartau$.
\end{theorem}

\begin{proof} The first statement follows from Theorem \ref{vertex}
and the fact that the origin is a vertex of $R_u$ if and only if $u$
is a vertex of $P_u$.  For the second claim we use Proposition
\ref{stdpair1}. The statement that $u$ is a vertex of $P_u$, and $u+v$
is a vertex of $P_{u+v}$ for all $v$ with $\mathrm{supp}(v) \subseteq
\tau$ is equivalent to the statement that the origin is a vertex of $R_u$
and it remains a vertex of $R_{u+v}$ for all such $v$.  Since
$\mathrm{supp}(v) \subseteq \tau$, this is the same thing as the
origin being the vertex of $R^{\bartau}_u$. Similarly, if for all $i
\not \in \tau$ there exists a $v'$ with $\mathrm{supp}(v') \subseteq
\tau \cup \{i\}$ such that the origin fails to be a vertex of
$P_{u+v'}$, then the origin is also not a vertex of $R_{u+v'}$, and hence
not a vertex of $R^{\bartau \setminus i}_u$, and vice versa.
\end{proof}

The characterization of the standard pairs in the above theorem 
also gives rise to a description of the irredundant irreducible 
primary decomposition of $\vrtx$. This is very similar to the 
description of the irredundant irreducible primary decompositions
of $\inI$ given in Section 4 of \cite{HT}. In order to give this
characterization we make the following definition.

\begin{definition} \label{critical}
 We call the polyhedron $Q_u$ {\em critical} if the origin is a vertex
of $R_u$, but not a
vertex of $R_{u+e_i}$ for any $i = 1, \ldots, k$.
\end{definition}

\begin{theorem} \label{primary}
The ideal $\vrtx$ has the irreducible primary decomposition
$$\vrtx = \bigcap_{Q^{\bartau}_u} \langle x_i^{u_i + 1}: \, 
i \in \bartau \rangle$$
where the intersection is taken over all critical $Q^{\bartau}_u$.
\end{theorem}

\begin{proof} The proof of Lemma 3.3 in \cite{STV} implies that
$$\vrtx = \bigcap_{(x^u, \tau)} \langle x_i^{u_i + 1}: \,
i \in \bartau \rangle$$
where the intersection is taken over all standard pairs $(x^u, \tau)$
such that $x^ux_i \in \vrtx$ for all $i \in \bartau$. By Theorem \ref{stdpair2}
these standard monomials are in bijection with critical
$Q^{\bartau}_u$.
\end{proof}

When we have a two-dimensional saturated lattice $\Lat = ker(A) \cap
\Int^n$, the codimension of $\vrtx$ is two, and so Proposition
\ref{dimension} implies that if an embedded prime $\mathcal{P}_{\tau}$
of $\vrtx$ exists, the codimension of $\mathcal{P}_{\tau}$ must be
three, which means $|\tau| = n-3$. Our next task is to show that in
this case $\mathrm{cone}\{a_i: i \in \tau\}$ cannot be a face of
$\mathrm{cone}\{a_i: i = 1,\ldots,n\}$ where $a_i$ is the $i$-th
column of the matrix $A$.  The result is a consequence of the
following lemma.

\begin{lemma} \label{polygon}
Let $Q \in \mathbf R^2$ be a polygon defined by $n$ facet-defining
inequalities $b_i \cdot x \leq u_i$, and let $R$ be the convex hull of
the lattice points in $Q$. Let $v$ be a vertex of $R$. Then there
exists a facet $j$ of $Q$ such that $v$ is a vertex of the convex
hull, $R_j$, of the lattice points in $ Q_j := \{x \in \mathbf R^2:
b_i \cdot x \leq u_i, \,\, i \neq j\}$.
\end{lemma}

\begin{figure}

\psfrag{v1}{$v_1$}
\psfrag{v2}{$v_2$}
\psfrag{negv1}{$-v_2$}
\psfrag{negv2}{$-v_1$}
\psfrag{v}{$v$}
\psfrag{k}{$K$}
\psfrag{negk}{$-K$}
\includegraphics{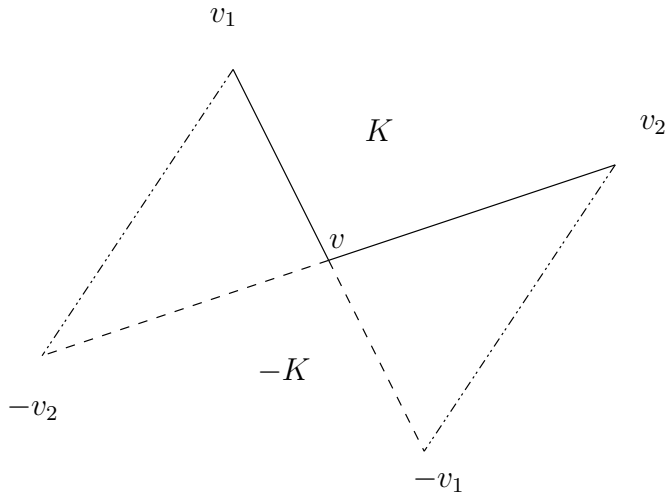}
\caption{The constructions of the proof of Lemma \ref{polygon} \label{proofeg}}
\end{figure}

\begin{proof} Suppose not. Clearly we can assume that $R$ is two-dimensional and that
$v$ is the origin. Let $v_1$ and $v_2$ be the two vertices of $R$
which are the neighboring vertices of the origin, in the clockwise and
counterclockwise directions respectively.  We define the pointed cone
$K$ generated by $v_1$ and $v_2$, and $-K$, the opposite cone
generated by $(-v_1)$ and $(-v_2)$.  These constructions are
illustrated in Figure \ref{proofeg}.  We first claim that each edge of
$Q$ has to intersect $-K$.  Suppose there is an edge $e$, lying on the
hyperplane $b_k \cdot x \leq u_k$, which does not intersect $-K$.
Then the convex region $S := \{x \in \mathbf R^2: b_i \cdot x \leq
u_i, \,\, i \neq k, \mathrm{and} \,\, b_k \cdot x \geq u_k\}$ does not
intersect $-K$ as well.  This is true because if $v \in S \cap -K$
there is a point $w$ on the line segment joining $v$ to the origin
lying on $e$, and $w$ would then be in $-K$.  Since the origin is not
a vertex of $\mathrm{conv}(Q_k \cap \mathbf Z^2)$, either $0$ is in
the interior of an edge of $R_k$ or it is in the interior of $R_k$. In
the first case there exists two vertices $y$ and $z$ of $R_k$ such
that $y \in S$ and $z \in R \subset K$ with $0 = \lambda y +
(1-\lambda)z$ for some $0 < \lambda < 1$. But then $y \in S \cap -K$,
contrary to our assumption. If $0$ is an interior point of $R_k$, then
there exist three vertices $y_1, y_2$ and $y_3$ of $R_k$ such that $0
= \lambda_1y_1 + \lambda_2y_2 + \lambda_3 y_3$ with $0 < \lambda_1,
\lambda_2, \lambda_3 < 1$ and $\sum \lambda_1 +\lambda_2 + \lambda_3
=1$.  Now, either exactly one or exactly two of these vertices are in
$S$.  In the first case, say $y_1 \in S$ and $y_2, y_3 \in R$, we have
$y_1 = -\frac{\lambda_2}{\lambda_1} y_2 - \frac{\lambda_3}{\lambda_1}
y_3$ and hence $y_1 \in (-K) \cap S$.  In the second case, say $y_1,
y_2 \in S$ and $y_3 \in R$, we have $ \frac{\lambda_1}{\lambda_1 +
\lambda_2} y_1 + \frac{\lambda_2}{\lambda_1 + \lambda_2} y_2 =
-\frac{\lambda_3}{\lambda_1 + \lambda_2}y_3$, and hence $
-\frac{\lambda_3}{\lambda_1 + \lambda_2}y_3 \in (-K) \cap S$.  In both
cases we get a contradiction to our assumption that edge $e$ does not
intersect $-K$.  This shows that all edges of $Q$ intersect $-K$.

Because $Q$ contains $v_1$ but not $-v_2$, and $v_2$ but not $-v_1$,
some edge of $Q$ must intersect the line segment $[v_1, -v_2]$, and
another one the line segment $[v_2, -v_1]$.  If we assume that the
facets of $Q$ are labeled going clockwise and the edge $1$ is the
first edge intersecting the facet of $-K$ defined by $(-v_1)$, then
edge 1 must be the edge intersecting $[v_2, -v_1]$. And if edge $n$ is
the last edge intersecting the facet of $-K$ defined by $(-v_2)$,
then edge $n$ must be the edge intersecting $[v_1, -v_2]$.  Edge 1 and
edge $n$ are the only edges of $Q$ not lying entirely in $-K$, so
they need to meet in a common vertex of $Q$.  But their endpoints
outside $-K$ are on opposite sides of the parallel line segments
$[v_1,-v_2]$ and $[v_2,-v_1]$, which makes this impossible.
\end{proof}

\begin{remark} Note that we cannot relax the hypothesis in Lemma \ref{polygon}
that $Q$ is a polygon to $Q$ being a possibly unbounded polyhedron.
An example of this phenomenon is in Figure \ref{unbounded}.  If any of
the facets of $Q$ are removed, the origin, $O$, ceases to be a vertex
of $R$.

\begin{figure}
\psfrag{O}{O}
\psfrag{Q}{Q}
\psfrag{R}{R}

\includegraphics[width=3in]{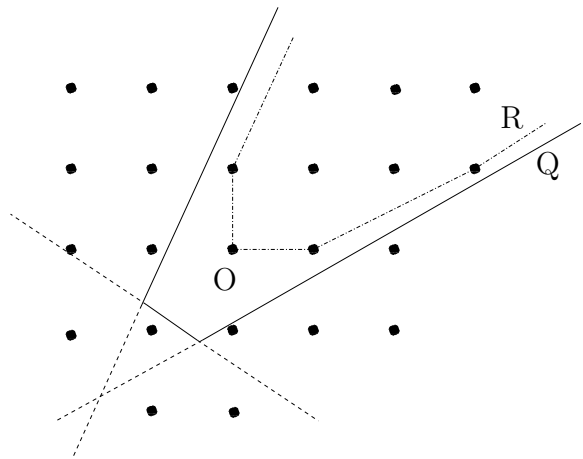}
\caption{A counterexample to Lemma \ref{polygon} for an unbounded polyhedron}
\label{unbounded}
\end{figure}

\end{remark} 

\begin{theorem}\label{codim2} Let $I_\Lat$ be a codimension two lattice 
ideal where $\Lat = ker(A) \cap \Int^n$ with $\Lat \cap \nint^n = \{0\}$. 
If $\mathcal{P}_\tau$ is an embedded
prime of $\vrtx$ then $\mathrm{cone}\{a_i: i \in \tau\}$ is not a face
of $\mathrm{cone}\{a_i: i = 1,\ldots,n\}$ where $a_i$ is the $i$-th
column of the matrix $A$.  In particular, the irrelevant maximal ideal
$\mathcal{P}_{\emptyset}$ is not associated to $\vrtx$.
\end{theorem}

\begin{proof} Let $(x^u, \tau)$ be a standard pair of $\vrtx$. Suppose that
$\mathrm{cone}\{a_i: i \in \tau\}$ is a face of $\mathrm{cone}\{a_i: i
= 1,\ldots,n\}$.  This means that the origin in $\mathbf R^2$ is in
the convex hull of $\{b_i : i \in \bartau \}$, where $b_i$ is the
$i$-th row of the $B$ defined after Proposition \ref{stdpair1}.  This
follows because positive covectors of (the oriented matroid of) $A$ 
correspond to positive
vectors of (the oriented matroid of) $B$ (see \cite[Chapter 6]{ziegler}).  
So $Q^{\bartau}_u$ is
a polygon. Theorem \ref{stdpair2} now implies that the origin in
$\mathbf R^2$ is a vertex of $R^{\bartau}_u$, but not a vertex of any
$R^{\bartau \setminus i}_u$ for $i \in \bartau$. But this is a contradiction
to Lemma \ref{polygon}.
\end{proof}

\begin{remark} The statement of Lemma \ref{polygon} also cannot be generalized to
higher dimensional polytopes. Similarly, the statement of Theorem
\ref{codim2} cannot be generalized to higher codimension. In
particular, the irrelevant maximal ideal $\mathcal{P}_{\emptyset}$
could be an embedded prime for some $\vrtx$ as the following example
shows.
\end{remark}

\begin{example} Let $A := [15,247, 248, 345]$. A lattice basis for $\Lat = ker(A) \cap 
\Int^4$ is given by the columns of 
$$B := \left[\begin{array}{rrr}
-4 & -3 & -3 \\
-6 & 9 & -2 \\
9 & -6 & -2 \\
-2 & -2 & 3 \\
\end{array}\right].$$
If we choose $u = (9,7, 7, 1)^T$, then $ Q^{\{1,2,3,4\}}_u = Q_u = \{x
\in \mathbf R^3: Bx \leq u\}$ is a tetrahedron. The polytope
$R^{\{1,2,3,4\}}_u = R_u$ has the following six vertices:
$$(0,0,-3) , (0,0,0) , (1,0,1) ,$$ 
$$(0,1,1) , (3,3,1),  (23,23,31).$$
Now $R^{\{1,3,4\}}_u$ contains the lattice point $(-1,0,-1)$ 
and the origin is in $\mathrm{conv}\{ (1,0,1) , (-1,0,-1)\}$; 
$R^{\{1,2,4\}}_u$ contains the lattice point $(0,-1,-1)$ and the origin is
in $\mathrm{conv}\{ (0,1,1) , 
(0,-1,-1)\}$;
$R^{\{1,2,3\}}_u$ contains the lattice points $(-1,0,0)$ and $(0,-1,0)$, and the origin is
in $\mathrm{conv}\{ (-1,0,0) , \\
(0,-1,0), (0,0,-1), (1,1,1)\}$, and finally
$R^{\{2,3,4\}}_u$ contains the lattice point $(-1,-1,-1)$ and the origin is
in $\mathrm{conv}\{ (1,1,1) , (-1,-1,-1)\}$. This shows the origin 
is not a vertex in any of these new polytopes. In particular,
$(x^u, \emptyset)$ is a standard pair of $\vrtx$, and hence the 
irrelevant ideal is an associated prime of $\vrtx$.
\end{example}

The above example also provides a counterexample to a conjecture about
the complexity of Gr\"obner fans of codimension three toric ideals
(Conjecture 6.2 in \cite{ST}). This conjecture stated that any
Gr\"obner cone of a codimension three toric ideal has at most four
facets.

\begin{theorem} \label{counterex}
There exists a toric ideal $I_A$ with $\mathrm{codim}(I_A) = 3$ which has 
a Gr\"obner cone with five facets.
\end{theorem} 

\begin{proof} Let $A$  be as in the above remark. 
If we choose $\omega = (111,0, 342, 1) $
as the cost vector we get the following reduced Gr\"obner basis:
$$
\begin{array}{l}
\{
a^{23}-d, 
da^{10}-bc, 
d^{12}a^4-b^{16}c, 
d^{55}a^3-b^{76}c,
d^{161}a^2-b^{225},\\

d^{204}a-b^{285},
d^{247}-b^{345}, 
cd^9a^7-b^{14},
cd^{20}a-b^{29},
cd^{63}-b^{89},\\

c^2d^8-b^{13}a^3,
c^4d^5-b^{11},
c^5d^4-b^{10}a^{10},
c^6d^2a^3-b^9, 
c^7a^{16}-b^8,\\

c^7d-b^8a^7,
c^8-b^7a^{17},
bca^{13}-d^2,
b^2c^2a^3-d^3,
b^3c^3-d^4a^7,\\

b^9a^{20}-c^6d^3,
b^{12}a^{13}-c^3d^7,
b^{15}a^6-d^{11},
b^{31}ca^2-d^{23},
b^{44}a^5-cd^{31},\\

b^{47}c^2-d^{35}a^2,
b^{60}a-d^{43},
b^{136}c-d^{98}a^2\}.
\end{array}$$
The corresponding Gr\"obner cone is given by 
$$\begin{array}{rrrrcr}
&+345b &  &  -247d   &  \leq & 0 \\    
 -20a &-9b & +6c & 3d    & \leq  & 0 \\    
 +2a & -136b & -c & +98d  & \leq & 0 \\   
 -3a & +76b & +c & -55d  & \leq & 0 \\   
 +7a & -3b & -3c & +4d  & \leq & 0 \end{array}$$
which are all facet defining. 
\end{proof}

This counterexample was found by using ${\tt TiGERS}$ \cite{HuT}, an
implementation to compute Gr\"obner fans of toric ideals developed by
Birkett Huber and Rekha Thomas.  Computer experiments with {\tt
TiGERS} have yielded many other examples of Gr\"obner cones of
codimension 3 toric ideals with five facets, and the following (thus far unique) codimension 3 toric ideal with a Gr\"obner cone with six facets.

\begin{example}
For the matrix
$$\left[ 
\begin{array}{ccccccc}
1 & 1 & 1 & 1 & 1 & 1 & 1 \\
2 & 8 & 9 & 7 & 10& 6 & 5 \\
8 & 7 & 4 & 8 & 7 & 2 & 2 \\
5 & 9 & 4 & 2 & 9 & 8 & 3\\
\end{array}
\right]
$$

the initial ideal of $I_A$ with respect to the weight vector $(252,
197, 0, 0, 153, 0, 0)$ corresponds to a Gr\"obner cone with six
facets.

\end{example}

\section{The Product Ideal}

In this section we define the {\em product ideal} of $\Lat$ which is
closely related to $\vrtx$, and which is much easier to compute.
Although in general the two ideals are not equal, we will look at two
special cases where they are: the case when $\Lat$ is unimodular, and
when $\Lat$ is a two-dimensional lattice in $\Int^2$.  Even in the
cases where they are not equal, we will show that the product ideal
carries valuable information about $\vrtx$. For instance we will show
that the radicals of the two ideals are equal.

\begin{definition} The product ideal $\prd$ is the monomial ideal
defined by 
$$\prd = \langle x^ux^v \, : \, u-v \in Gr_\Lat \rangle.$$
\end{definition}

Since each initial ideal $\inI$ contains one of $x^u$ or $x^v$
whenever $u-v \in \Lat$, we have $\prd \subseteq \vrtx$.  This containment
could be strict, however, as shown in the following example.  Let
$\Lat = ker(A) \cap \Int^3$ where $A = [3 \, \, 4 \, \,5]$. Then $\prd
= \langle ab^2c, a^2bc^2, a^3bc, a^4b^3, a^5c^3, b^5c^4 \rangle$ is
strictly contained in $\vrtx = \langle ab^2c, a^2bc, a^4b^3, a^5c^3,
b^5c^4 \rangle$. There are two special cases, though, in which the
product ideal and the vertex ideal are equal.  The first case is when
$\Lat$ comes from a unimodular matrix.  We recall that a $d \times n$
matrix is unimodular if all maximal $d \times d$ minors have the same
absolute value.

\begin{proposition}\label{uni}
If $\Lat = ker(A) \cap \Int^n$ where $A$ is a unimodular matrix, then
$\prd = \vrtx$, and $\prd$ coincides with the matroid ideal 
$I_{\Delta(\mathcal{M}(\Lat))}$.
\end{proposition}

\begin{proof}The initial ideals $\inI$ are all square-free (Corollary 8.9 in
\cite{GBCP}) and hence $\vrtx$ is radical.  Therefore, by Corollary
\ref{matroid}, $\vrtx = I_{\Delta(\mathcal{M}(\Lat))}$. 
But the minimal generators of $I_{\Delta(\mathcal{M}(\Lat))}$ are of the
form $\prod_{i \in J} x_i$ for some $J=\{i_1,\dots, i_k \}$ such that
$\{a_{i_1}, \ldots, a_{i_k}\}$ is a circuit of $A$. Now Proposition
8.11 of \cite{GBCP} implies that the Graver basis of $I_\Lat$ is
$\{\alpha_i-\beta_i \in \Lat\, : \, \mathrm{supp}(\alpha_i-\beta_i)$
is the support of a circuit $\}$.  Therefore $\prd = \langle\prod_{i \in J}
x_i : J=\{i_1,\ldots, i_k\}$ is the support of a circuit $\rangle$, and
hence $\prd = \vrtx = I_{\Delta(\mathcal{M}(\Lat))}$.
\end{proof}

\begin{proposition}\label{dim2} If $\Lat$ is a two-dimensional lattice 
in $\Int^2$, then $\prd = \vrtx$. 
\end{proposition}

\begin{proof} Let $S = k[x,y]$ and suppose $x^uy^v \in \vrtx$, so  
$(u,v)$ is not a vertex of its fiber $P_{(u,v)}$, but $x^uy^v \not \in
\prd$.  If $(a,b) \in P_{(u,v)}$ where $(0,0) \leq (a,b) \leq (u,v)$,
then $x^uy^v \in \prd$ because $x^{u-a}y^{v-b} -1 \in I_\Lat$ and
hence $x^{u-a}y^{v-b} \in \prd$.  So no such point in $P_{(u,v)}$
exists.  Now there must be a vertex $(a,b)$ of this fiber with $b <
v$, because otherwise $(u,v)$ would be a vertex. Let $(a,b)$ be the
vertex with $b < v$ such that $(a,b)$ is the maximum with this
property. Let $H$ be the line through $(u,v)$ and $(a,b)$, let $H^-$
be the halfspace containing the origin, and let $H^+$ be the other
halfspace. If $P_{(u,v)} \subseteq H^+$, since $(u,v)$ is not a vertex
of $P_{(u,v)}$, the line $H$ must contain $(c,d) \in P_{(u,v)}$ such
that $0 \leq c \leq u$ and $v < d \leq 2v-b$.  But then $(u-c, v-d)
\in \Lat$ and $x^{u-c}-y^{d-v} \in I_\Lat$, which implies
$x^{u-c}y^{d-v} \in \prd$. This implies $x^uy^v \in \prd$ since $d-v
\leq v$. Hence we are reduced to the case that $P_{(u,v)}$ is not
contained in $H^+$ and $2u<a$ (so no such $(c,d) \in P_{(u,v)})$.
This means that there exists a vertex $(e,f) \in H^- \cap P_{(u,v)}$.
Now if $e>u$, by the construction of $(a,b)$ we must have $f<b$.  If
in addition $e<a$, the existence of a vector of the form $(k,0) \in
\Lat$ for some $k$ means that $(e+kN,f) \in P_{(u,v)}$ for $N \gg 0$
which contradicts $(a,b)$ being a vertex.  On the other hand, if
$a<e$, $(a,b)$ would not be a vertex of $P_{(u,v)}$.  So we conclude
that $e<u$.  But now we must have $v<f<2v-b$, where the second
inequality follows from the assumption that $(e,f) \in H^-$ and
$2u<a$.  Since $(e-u,f-v) \in \Lat$, it follows that $x^uy^v \in
\prd$, a contradiction, so $\prd=\vrtx$.
\end{proof}

\begin{example} \label{dim3ex}
Proposition \ref{dim2} fails when $dim(\Lat) \geq 3$. For instance,
let $\Lat$ be the lattice in $\Int^3$ generated by the columns of the
matrix
$$\left[\begin{array}{ccc}
1 & 4 & 3 \\
-2 & 0 & 5 \\
-1 & 1 & -9 \end{array}\right].$$

One can verify using {\tt Macaulay2} \cite{M2} that $\prd = \langle
ab^2c, a^4c, a^5b^2, b^8c^5, \\ abc^{12}, b^3c^{11}, b^{19}c, ab^{21},
a^4b^{19}, ac^{26}, a^3c^{25}, b^2c^{27}, bc^{38}, a^{49}b, c^{103},
b^{103}, a^{103} \rangle$, and it is strictly contained in $ \vrtx =
\langle c^3, ab^2c, a^4c, a^5b^2, b^{19}c, ab^{21}, a^4b^{19},
a^{49}b,\\ 
b^{103}, a^{103} \rangle$.
\end{example}
The example at the beginning of this section shows that Proposition \ref{dim2}
does not hold even for a two-dimensional lattice $\Lat$  when $\Lat$ is in
$\Int^n$ for $n \geq 3$. However, we will show that $\prd$ and $\vrtx$ have
the same radical, and that for two dimensional lattices they are almost
equal.

For $\sigma \subseteq \{1,\ldots,n\}$, let $\pi_\sigma: \Int^n \longrightarrow
\Int^{n-|\sigma|}$ be the projection map which eliminates the coordinates
indexed by $\sigma$. We will denote the image of a lattice $\Lat$ under
this map by $\Lat_\sigma$. It is clear that if $\mathrm{dim}(\Lat) = \mathrm{dim}(\Lat_\sigma)$
then $\Lat$ and $\Lat_\sigma$ are isomorphic lattices. This observation 
implies the following lemma.

\begin{lemma} \label{local-Graver}
Let $Gr_\Lat$ and $Gr_{\Lat_\sigma}$ be the Graver bases of the lattices
$\Lat$ and $\Lat_\sigma$ where $\mathrm{dim}(\Lat) = \mathrm{dim}(\Lat_\sigma)$. Then
$Gr_{\Lat_\sigma} \subseteq \pi_\sigma(Gr_\Lat)$. 
\end{lemma}

\begin{proof}
If $\alpha^{\prime}-\beta^{\prime} \in Gr_{\Lat_\sigma}$, there is a
unique $\alpha-\beta \in \Lat$ such that $\pi_\sigma(\alpha-\beta) =
\alpha^{\prime}-\beta^{\prime}$. If $\alpha=\alpha_1+\alpha_2$,
$\beta=\beta_1+\beta_2$, where $\alpha_i-\beta_i \in \Lat$ for
$i=1,2$, and $\alpha_i,\beta_i \in \mathbf N^n$, then
$\alpha^{\prime}=\pi_{\sigma}(\alpha_1)+\pi_{\sigma}(\alpha_2)$ and
$\beta^{\prime}=\pi_{\sigma}(\beta_1)+\pi_{\sigma}(\beta_2)$, with
$\pi_{\sigma}(\alpha_i)-\pi_{\sigma}(\beta_i) \in \Lat_{\sigma}$ for
$i=1,2$.  As this contradicts $\alpha^{\prime}-\beta^{\prime} \in
Gr_{\Lat_\sigma}$, we conclude that $\alpha-\beta \in Gr_{\Lat}$, so
$\alpha^{\prime} -\beta^{\prime} \in \pi_{\sigma}(Gr_{\Lat})$.
\end{proof}

The algebraic analogue of the projection map $\pi_\sigma$ is the
localization map $ \hat{\pi}_\sigma : k[x_1,\ldots, x_n]
\longrightarrow k[x_i: \, i \notin \sigma]$ where
$\hat{\pi}_\sigma(x_i) = x_i$ if $i \notin \sigma$ and
$\hat{\pi}_\sigma(x_i) = 1$ otherwise. This corresponds to localizing
at the monomial prime $\mathcal{P}_\sigma = \langle x_i \, : \, i \notin \sigma
\rangle$.  We now compare $Top(\prd)$ with $Top(\vrtx)$, where
$Top(M)$ is the intersection of the top-dimensional primary components
of the ideal $M$. When we consider a monomial ideal $M$ with
top-dimensional minimal primes $\mathcal{P}_{\sigma_1}, \ldots, \mathcal{P}_{\sigma_k}$,
 we have $Top(M) = \cap_{i=1}^{k} \hat{\pi}_{\sigma_i}(M)$.

\begin{proposition} \label{localize} 
If $\mathrm{dim}(\Lat)=\mathrm{dim}(\Lat_\sigma)$ then $\hat{\pi}_\sigma(\vrtx) = V_{\Lat_\sigma}$ and 
$\hat{\pi}_\sigma(\prd) = P_{\Lat_\sigma}$. 
\end{proposition}

\begin{proof}
From Lemma \ref{local-Graver} we know that $\pi_\sigma(Gr_\Lat) 
\subseteq Gr_{\Lat_\sigma}$.
Let $x^u$ be a minimal generator of $\vrtx$.  By Corollary
\ref{mingens} we know that $x^u=\lcm_i(x^{\alpha_i})$ where $\sum_i c_i
(\alpha_i-\beta_i)=0$ for $\alpha_i-\beta_i \in Gr_{\Lat}$ and $c_i > 0$.
Now $\sum_i c_i
\pi_{\sigma}(\alpha_i-\beta_i)=0$.  Writing
$\pi_{\sigma}(\alpha_i-\beta_i)=\sum_j(\alpha_{ij}-\beta_{ij})$ where
$\alpha_{ij}-\beta_{ij} \in Gr_{\Lat_{\sigma}}$ and $\alpha_{ij} \leq
\alpha_i$, $\beta_{ij} \leq \beta_i$ for all $j$, we see that
for $x^v=\lcm_{i,j}(x^{\alpha_{ij}}), \, x^v \in V_{\Lat_{\sigma}}$.  
Since $x^v$ divides $ \lcm_i(\alpha_i)$, it follows that
$\hat{\pi}_{\sigma}(x^u) \in V_{\Lat_{\sigma}}$.

For the other inclusion, let $x^u$ be a minimal generator of
$V_{\Lat_{\sigma}}$, so $x^u=\lcm_i(x^{\alpha_i})$ for $\sum_i c_i (
\alpha_i-\beta_i)=0$, where $\alpha_i-\beta_i \in Gr_{\Lat_{\sigma}}$
and $c_i > 0$.  Let $\alpha_i^{\prime} - \beta_i^{\prime}$ be the preimage
of $\alpha_i - \beta_i$ under $\pi_{\sigma}$.  We still have $\sum_i c_i (
\alpha_i^{\prime}-\beta_i^{\prime})=0$, so for
$x^v=\lcm_i(x^{\alpha_i^{\prime}})$, $x^v \in \vrtx$, and thus
$\hat{\pi}_{\sigma}(x^v)=x^u \in \hat{\pi}_{\sigma}(\vrtx)$.

The second statement of the proposition follows from the definition of
the product ideal, and the observation that if
$\pi_{\sigma}(\alpha-\beta) \not \in Gr_{\Lat_{\sigma}}$ for
$\alpha-\beta \in Gr_{\Lat}$, we can write $\pi_{\sigma}(\alpha-\beta)$
as the sum of $\alpha_i-\beta_i \in Gr_{\Lat_{\sigma}}$ so that
$x^{\alpha_i+\beta_i} |
x^{\pi_{\sigma}(\alpha)+\pi_{\sigma}(\beta)}=\hat{\pi}_{\sigma}(x^{\alpha+\beta})$.
\end{proof}

\begin{cor} \label{rad-equal}
The radical of $\prd$ and the radical of $\vrtx$ coincide.  Moreover,
$Top(\prd) \subseteq Top(\vrtx)$.
\end{cor}

\begin{proof} 

Theorem \ref{radva} shows that $rad(\vrtx)$ is an equidimensional
ideal.  Now an associated prime $\mathcal{P}_\sigma$ of $\vrtx$ is a minimal
prime if and only if $\Lat_\sigma$ is a full dimensional lattice in
$\mathbf Z^{n-|\sigma|}$.  But this is true if and only if there exist
$n_i$ such that $n_ie_i \in \Lat_{\sigma}$ for all $i \not \in
\sigma$.  This happens if and only if $x_i^{n_i} \in P_{\Lat_\sigma}$
for all $i \not \in \sigma$, which happens exactly whenever
$P_{\Lat_\sigma} = \hat{\pi}_\sigma(P_\Lat)$ is a zero-dimensional
ideal, and hence $\mathcal{P}_\sigma$ is a minimal prime of $P_\Lat$.  This
shows that $rad(\vrtx) = rad(\prd)$.

The second statement follows from Proposition \ref{localize} and the
 discussion before it, and the fact that $P_{\Lat_\sigma} \subseteq
 V_{\Lat_\sigma}$.
\end{proof}

\begin{cor}\label{codimtwo}
If $\mathrm{dim}(\Lat) = 2$ then $Top(\prd) = Top(\vrtx)$.
\end{cor}

\begin{proof} 
Proposition \ref{dim2} says that $P_{\Lat_\sigma} = V_{\Lat_\sigma}$
when $\mathcal{P}_\sigma$ is a minimal prime of $\vrtx$ (and of
$\prd$).  Now Proposition \ref{localize}, Corollary \ref{rad-equal}
and the discussion before them imply the result.
\end{proof}
We note that the above corollary fails when $\mathrm{dim}(\Lat) \geq 3$.
Example \ref{dim3ex} provides a lattice $\Lat \in \Int^3$ of dimension
three. Therefore $Top(\vrtx) = \vrtx$ and $Top(\prd) = \prd$, but in that
example we saw that $\prd \neq \vrtx$.

\end{document}